\documentclass[english]{amsart}

\usepackage{babel}
\usepackage{amstext}
\usepackage{amsmath}
\usepackage{amsfonts}
\usepackage{latexsym}
\usepackage{ifthen}
\usepackage{xypic}
\xyoption{all}
\newcommand\sG{{\mathcal G}}

\newcommand\sO{{\mathcal O}}

\newcommand\sC{{\mathcal C}}

\newcommand\Aut{{\rm Aut}}
\newcommand\Hol{{\rm Hol}}

\newcounter{lemma}

\newtheorem{lemma1}[lemma]{\setcounter{equation}{0}}

\newenvironment{lemma}{\begin{lemma1}{\bf Lemma.}}{\end{lemma1}}

\newenvironment{theorem}{\begin{lemma1}{\bf Theorem.}}{\end{lemma1}}

\newenvironment{proposition}{\begin{lemma1}{\bf Proposition.}}{\end{lemma1}}

\newenvironment{construction}{\begin{lemma1}{\bf Construction.}}{\end{lemma1}}

\newenvironment{Induction Step}{\begin{lemma1}{\bf Induction Step.}}{\end{lemma1}}
\newenvironment{Proof of Theorem 1.2}{\begin{lemma1}{\bf Proof of Theorem 1.2.}}{\end{lemma1}}

\makeatletter
\ifnum\@ptsize=0  \fi \ifnum\@ptsize=2  \fi 

\title {Holomorphic maps onto K\"ahler manifolds with non-negative Kodaira dimension} \author{Jun-Muk Hwang and  Thomas Peternell}
\date{today}
\date{\today, {\em Address:}   Korea Institute for Advanced Study, 207-43, 
Cheongnyangni-dong, Seoul, 130-722, Korea,
{\em email:}jmhwang@kias.re.kr; Mathematisches Institut, Universit\"at Bayreuth, D--95440 Bayreuth, Germany, {\em e-mail:} thomas.peternell@uni-bayreuth.de}

\begin{document}

\begin{abstract}
This paper studies the deformation theory of a holomorphic
surjective map from a normal compact complex space $X$ to a
compact K\"ahler manifold $Y$. We will show that when the target
has non-negative Kodaira dimension,  all deformations of
surjective holomorphic maps $X \to Y$ are
  unobstructed, and the associated components of $\Hol(X,Y)$ are complex
  tori. Under the additional assumption that  $Y$ is  projective algebraic, this was
  proved  in [HKP06]. The proof in [HKP06] uses the algebraicity
  in an essential way and cannot be generalized directly to the K\"ahler
  setting. A new ingredient here is a careful study
of  the infinitesimal deformation of  orbits of an action of a
complex torus. This study, combined with the result for the
algebraic case, gives the proof for the K\"ahler setting.
\end{abstract}

\maketitle

\tableofcontents

\section{Introduction}

This paper studies the deformation theory of a holomorphic
surjective map from a normal compact complex space to a compact
K\"ahler manifold. For a holomorphic map $f:X \to Y$ between two
compact complex spaces, denote by $\Hol(X,Y)$ the space of
holomorphic maps from $X$ to $Y$ and by $\Hol_f(X,Y)$ the
connected component of $\Hol(X, Y)$ containing the point
corresponding to $f$. We prove:

\begin{theorem} Let $f: X \to Y$ be a surjective map between a normal compact complex space $X$ and a compact K\"ahler manifold $Y$
of nonnegative Kodaira dimension. Then there exists a
factorization
 $$
  \xymatrix{
    X  \ar[r]_{\alpha} \ar@/^0.3cm/[rr]^{f} &
    Z  \ar[r]_{\beta} &
    Y
  }
  $$
  where
  \begin{enumerate}
  \item $\beta$ is a finite unramified covering and
\item if $\Aut^0(Z)$ is the maximal connected subgroup of the
    automorphism group of $Z$, then $\Aut^0(Z)$ is a complex torus,
    and the natural morphism
    $$
     \Aut^0(Z) \Big / \Aut(Z/Y) \cap \Aut^0(Z)  \to \Hol_f(X,Y)
    $$
    is isomorphic.
  \end{enumerate}
  In particular, all deformations of surjective holomorphic maps $X \to Y$ are
  unobstructed, and the associated components of $\Hol(X,Y)$ are complex tori.
\end{theorem}

\vskip .2cm \noindent This answers the question raised in [Hw05,
p. 769] and [HKP06, Remark 1.4] affirmatively.

\vskip .2cm
For simplicity, let us say that a surjective
holomorphic map $f:X \to Y$ from a normal compact complex space
$X$ to a compact complex manifold $Y$ is {\it rigid} if there exists $Z$
factorizing $f$ as in (1.1) such that
$$ g^* H^0(Z, h^*T_Y) =
g^* H^0( Z, T_Z).$$
   Since the tangent space to $\Hol(X,Y)$
at $f$ is $H^0(X,f^*T_Y),$  we can reformulate Theorem 1.1:

\begin{theorem} Let $X$ be a normal compact complex space and $Y$ be
a compact K\"ahler manifold of non-negative Kodaira dimension.
Then $f$ is rigid.
\end{theorem}

Note that it suffices to prove Theorem 1.2 for an individual $\sigma \in H^0(X,f^*T_Y):$
there exists an unramified holomorphic covering
$h: Z \rightarrow Y$ with  a holomorphic map $g:X \rightarrow Z$
satisfying
$f= h \circ g$, such that $$ \sigma \in  g^* H^0(Z, h^*T_Y) =
g^* H^0( Z, T_Z).$$ Here of course $Z$ and $h$ depend on $\sigma,$ but it is clear by taking fiber products that
one can choose $Z$ and $h$ independent from $\sigma.$

In case $X$ and $Y$ are projective, Theorem 1.1, and equivalently,
Theorem 1.2,  was proved in [HKP06]. The proof of Theorem 1.2 in
the projective setting in [HKP06] depends on Miyaoka's
semi-positivity theorem whose proof in turn requires the use of
characteristic $p>0$ method. So the proof in [HKP06] cannot be
generalized directly to compact K\"ahler
manifolds. \\
Theorem 1.2 was proved for a compact K\"ahler manifold $Y$ with
trivial canonical class in [Hw05]. The proof in [Hw05], which
generalizes the earlier work of [KSW81],  uses the differential
geometry of a Ricci-flat metric. The proof in [Hw05] seems
difficult to generalize to prove Theorem 1.2, because the
non-negativity of Kodaira dimension alone is too weak to give a
nice K\"ahler metric.

The key idea of this paper is that the following weaker
version of Theorem 1.2 can be proved without characteristic $p>0$
or differential geometric techniques.

\begin{theorem}
Let $X$ be a normal compact complex space and $Y$ be
a compact K\"ahler manifold of non-negative Kodaira dimension.
Given a surjective holomorphic map  $f: X
\rightarrow Y$ and  a section $$\sigma \in df(H^0(X, T(X))) \subset
 H^0(X, f^*T(Y)),$$
  there exists an unramified holomorphic covering
$h: Z \rightarrow Y$ with  a holomorphic map $g:X \rightarrow Z$
satisfying
$f= h \circ g$, such that $$ \sigma \in  g^* H^0(Z, h^*T(Y)) =
g^* H^0( Z, T(Z)).$$
\end{theorem}

\medskip
Theorem 1.3 will be proved in Section 3. The main idea is to use
the infinitesimal deformation of holomorphic maps from a complex
torus to $Y$ arising from the section $\sigma$ and $f$. The
unramified covering is constructed by showing that an essential
part of
 the space of these holomorphic maps is  smooth.

\medskip
Once Theorem 1.3 is established, we will show that Theorem 1.2 in the projective case and Theorem
1.3
implies Theorem 1.2 (Section 4). In this sense, the proof of Theorem 1.1
does depend on the method of characteristic $p>0$.

\medskip The authors were supported by the DFG priority program ``Global methods in Complex Geometry'', which we gratefully
acknowledge.

\section{Preliminaries}
\setcounter{lemma}{0}

\begin{lemma} Let $Y$ be a compact K\"ahler manifold with a complex torus $G$ acting effectively on $Y$ as a group of holomorphic
transformations. Then there exists a compact K\"ahler manifold
$Y/G$ and a submersion $q: Y \to Y/G$ such that the fibers of $q$
are the $G-$orbits. In particular, $Y$ and $Y/G$ have the same
Kodaira dimension.
\end{lemma}

The proof is found in [Cr72, p.225]. We call $q: Y \to Y/G$ the
$G-$quotient of $Y.$

\begin{lemma} In the situation of Lemma 2.1 let $Z = Y/G.$ Suppose furthermore that there is an unramified covering $s: \tilde Z \to Z$ and
set $\tilde Y = Y \times_Z \tilde Z$ with projection $\tilde q:
\tilde Y \to \tilde Z.$ Then we have an exact sequence
$$ 0 \to H^0(\tilde Y,T_{\tilde Y/\tilde Z}) \to H^0(\tilde Y,T_{\tilde Y}) \to H^0(\tilde Y,\tilde q^*T_{\tilde Z}) \to 0.$$
\end{lemma}

\begin{proof} It suffices to notice that the connecting homomorphism, i.e., the Kodaira-Spencer map,
$$ H^0(Y,\tilde q^*T_{\tilde Z}) \to H^1(\tilde Y,T_{\tilde Y/ \tilde Z}) $$
vanishes, since fibers of $\tilde q$ are $G-$orbits.
\end{proof}

\begin{proposition} Let $q: Y \to Z$ be a $G-$quotient of $Y$ as in Lemma 2.1. Then a surjective holomorphic map $f: X \to Y$ is
rigid if $q \circ f: X \to Y/G$ is rigid.
\end{proposition}

\begin{proof} Since $p := q \circ f$ is rigid,
there is an unramified covering $s: \tilde Z \to Z$ and a factorization
 $$
  \xymatrix{
    X  \ar[r]_{r} \ar@/^0.3cm/[rr]^{p} &
    \tilde Z \ar[r]_{s} &
    Z
  }
  $$ such that
$$H^0(X,p^*T_Z) = r^*H^0(T_{\tilde Z}).$$
 Let $\tilde Y = Y \times_Z \tilde Z$ with projection $\tilde q: \tilde Y \to \tilde Z.$ Let $\tilde f: X \to \tilde Y$ be the
canonical map. \\
Notice that $T_{Y/Z} $ is spanned by the $G-$vector fields and
therefore is trivial. Hence $T_{\tilde Y/\tilde Z}$ is trivial as
well, which shows
$$ H^0(X,\tilde f^*T_{\tilde Y/\tilde Z}) = \tilde f^*H^0(\tilde Y,T_{\tilde Y/\tilde Z}).$$
Now a diagram chase in the following diagram yields the claim
(the surjectivity in the first row comes from Lemma 2.2).
$$\begin{array}{lllllcl} H^0(\tilde Y,T_{\tilde Y/\tilde Z}) &\hookrightarrow & H^0(\tilde Y,T_{\tilde Y}) & \to &H^0(\tilde Y,\tilde q^*T_{\tilde Z}) & \to & 0\\
\downarrow & & \downarrow & & \downarrow & & \\
H^0(X,\tilde f^T_{\tilde Y/\tilde Z}) & \hookrightarrow &H^0(X,\tilde f^*T_{\tilde Y}) & \to & H^0(X,\tilde f^*\tilde q^*T_{\tilde Z}) & \to &
H^1(X,\tilde f^*T_{\tilde Y/\tilde Z}) \\
\end{array} $$
\end{proof}

The following proposition is obvious from the definition of
rigidity.

\begin{proposition}
A surjective holomorphic map  $f: X \to Y$ from a normal compact
complex space to a compact complex manifold is rigid if there
exists a factorization
 $$
  \xymatrix{
    X  \ar[r]_{\alpha} \ar@/^0.3cm/[rr]^{f} &
    Z  \ar[r]_{\beta} &
    Y
  }
  $$
with $\beta $ unramified and $\alpha$ rigid.
\end{proposition}

We will need the following two results of C.Horst. Proposition 2.5
is [Ho87, Theorem 0.2.1] and Proposition 2.6 is [Ho87, Corollary
5.1.1]. Note that in these two propositions, the factorizing map
$h$ is {\it not} required to be unramified, which is the essential
difference between [Ho87] and our Theorem 1.3.

\begin{proposition} Let $f: X \to Y$ be a finite surjective map
 between a normal compact complex space $X$ and a compact complex manifold $Y$ of nonnegative Kordaira dimension.
 Let $Z \subset \Hol_s(X,Y)$
be a compact subvariety of the space of surjective holomorphic maps. 
Then there exists a
factorization
 $$
  \xymatrix{
    X  \ar[r]_{k} \ar@/^0.3cm/[rr]^{f} &
    Y'  \ar[r]_{h} &
    Y
  }
  $$
such that
$$ Z \subset  h  \circ \Aut^0(Y') \circ k.$$
\end{proposition}

The second proposition replaces the compactness assumption by a condition on the branch locus. 
\begin{proposition} Let $f: X \to Y$ be a finite surjective map between a normal compact complex space $X$ and a compact complex manifold $Y$.
 Let $Z \subset \Hol_f(X,Y)$
be an irreducible closed subvariety containing $[f]$. For each $g
\in Z$, denote by $B_g \subset Y$ the branch locus of $g$. Suppose
that $B_g$ and $g^{-1}(B_g)$ are independent of $g$ for all $g \in
Z$. Then there exists a factorization
 $$
  \xymatrix{
    X  \ar[r]_{k} \ar@/^0.3cm/[rr]^{f} &
    Y'  \ar[r]_{h} &
    Y
  }
  $$
such that
$$ Z \subset  h  \circ \Aut^0(Y') \circ k.$$
\end{proposition}

\section{Proof of Theorem 1.3}
\setcounter{lemma}{0}

In this section we prove Theorem 1.3. So we consider the differential
$$ df: T_X \to f^*T_Y $$
and the associated map
$$ df: H^0(X,T_X) \to H^0(Y,f^*T_Y).$$
We consider
$$ 0 \neq v \in H^0(f^*T_Y) $$
and assume that $v$ is of the form
$$ v = df(v') $$
with $$v' \in H^0(X,T_X).$$
By integration the automorphism group of $X$ is positive-dimensional, and since $X$ is not uniruled, the identity component is a torus.
First we show

\begin{proposition}
Let $X$ be a normal compact  complex variety and $G$ be a complex torus acting
on $X$ effectively by a holomorphic map $$\Phi: X \times G
\rightarrow X.$$
 Let $Y$ be a compact complex manifold and $f: X \rightarrow Y$ be a
finite surjective holomorphic map. Denote by $F: X \times G
\rightarrow Y$ the composition of
$\Phi$ and $f$. For $x \in X$, let $$F_x: G \rightarrow Y$$ be the map
defined by $F_x(g) := F(x, g)$ for $g \in G$.  Assume that there exists a
subvariety $E$ of codimension $\geq 2$ in $X$ which is preserved by
$G$, namely, $\Phi(E \times G) = E$, such that for each $x
\in X \setminus E$,  the pull-back $F_x^*T_Y$
 is a trivial
vector bundle on $G$. \\
Then there exists an
unramified covering $ Y_1 \rightarrow Y$ factoring $f$ such that the
$G$-action descends to $Y_1$.
\end{proposition}

\begin{proof}
Note that $F(E \times G) = f(E)$ is of codimension $\geq 2$ in $Y$.
We have a natural holomorphic
map
$$\chi: X \rightarrow \Hol(G, Y)$$
defined by $\chi(x) = F_x$.
 Let $B \subset \Hol(G, Y)$ be the image
 $\chi(X)$. For the identity element $e \in G$, the map
$\Phi(\cdot, e): X \rightarrow X$ is the identity map of $X$. Thus
when $b = \chi(x)$, we have $$b(e) = f \circ \Phi(x, e) = f(x)$$
and  $$\{ b(e) \in Y, b \in B \} = \{f(x), x \in X \} = Y.$$
 This implies that  the irreducible variety $B$ has dimension $n = \dim Y$.
\vskip .2cm \noindent
For $x \in X \setminus E$, let
$b:= \chi(x): G \rightarrow Y.$ Then  $b^* T_Y= F_x^*T_Y$ which
 is a trivial bundle  by the
assumption.  It follows that the tangent
space  (e.g. [Ko96, I.2.16])
to $\Hol(G, Y)$ at $b$,  $H^0(G, b^*T_Y)= H^0(G, \sO_G^n)$, has dimension $n$.
 Thus $b$ is a smooth point of $B$ and $B$ must be an irreducible
  component of $\Hol(G, Y)$.
\vskip .2cm \noindent
Now let $\Psi: B \times G \rightarrow Y$ be
 the restriction of the evaluation map
$$\Hol(G, Y) \times G \rightarrow Y.$$ Given $ b\in B$ and $p \in G$, let
$$\phi(p,b): T_b(B) = H^0(G, b^*T(Y)) \rightarrow b^*T_pY$$ be the evaluation
at $p$. As in [Ko96, II.3.4],  the differential $$d\Psi_{b,p}:
T_b(B) \times T_p(G) \rightarrow T_yY, \; y:= \Psi(b,p),$$ is $
\phi(p,b) + db(p)$. Suppose $b = \chi(x), x \in X \setminus E$, By
the assumption on the triviality of $b^*T(Y)$, we see that
$\phi(p,b)$ is surjective. It follows that $\Psi$ is submersive
 at $b$. (This part is analogous to [Ko96, II.3.5]).
In particular, the fiber $\Psi^{-1}(y)$ is smooth when $y \in Y \setminus f(E)$.
Now let $j: B\times G \rightarrow Y_1$ and $h: Y_1 \rightarrow Y$ be
the Stein factorization of $\Psi$. Then $h$
must be unramified over $Y \setminus f(E)$, hence unramified
everywhere.
 Since $j$ has connected fibers, the natural
$G$-action on $B \times G$ must descend to a $G$-action on $Y_1$.
Composing $$X \times \{ e \}  \subset X \times G \rightarrow B
\times G \rightarrow Y_1 \rightarrow Y,$$ we get a factorization
of $f:X \rightarrow Y$ with the desired properties.
\end{proof}

We have to verify that the technical assumption in Proposition 3.1 is indeed true:

\begin{proposition}  Let $X$ be a normal compact complex variety and $G$ be a
 complex torus acting on $X$. Let $Y$ be a
 compact
complex manifold with non-negative Kodaira dimension and  $f: X
 \rightarrow Y$ be a  finite holomorphic map. Let  $\Phi, F, F_x$
be as in the previous proposition. Then there exists a subvariety
$E \subset X$ of codimension $\geq 2$ which is preserved by $G$, such
that for any
 $x \in X \setminus E$,  $F_x^*T_Y$
 is a trivial vector bundle on $G$.
\end{proposition}

\begin{proof}  Since $X$ is normal, we may put all singular points of $X$ in
our $E$. Thus we may consider only smooth points of $X$. Let $x\in X$
be a smooth point. Choose a small neighborhood $U$ of $x$
where $T_X$ is a trivial bundle. Choose a holomorphic frame
$v_1, \ldots, v_n$ of $T_X|_U$ and regard them as vector fields
on $U \times G$ by the projection $U \times G \rightarrow U$. Let
$$u_i := dF(v_i)  \in H^0(U \times G, F^*T_Y), i= 1, \ldots, n.$$
We will use these sections of $F^*T_Y$ on $U \times G$ to show
that $F_x^*T_Y$ is trivial for suitable choices of $x$. \vskip
.2cm \noindent First, we will show this for any smooth  point $x
\in X$ such that $\Phi(x,g) = g\cdot x \not\in Ram(f)$ for some $g
\in G$, where $Ram(f)$ denotes
 the underlying reduced divisor of the ramification of
$f$. In the construction of
 the sections $u_1, \ldots, u_n$, we may assume that $g \cdot U$ is
disjoint from $Ram(f)$ by shrinking $U$. Then $u_1, \ldots, u_{n}$
 will be pointwise
independent at every point of $ U \times \{g\}$. This implies that
$u_1 \wedge \cdots \wedge u_n$ defines a section of $F^*K_Y^{-1}$
on $U \times G$ whose zero divisor $Z$ is disjoint from $U \times
\{g\}$. Since the Kodaira dimension of $Y$ is non-negative,
$F_z^*K_Y^{-1}$ cannot have a non-zero section with non-empty zero
for general $z \in U$. Thus  $Z$ is disjoint from $\{z\} \times G$
for general $z \in U$. It follows that $Z$ is empty and $u_1,
\ldots, u_n$ are pointwise independent everywhere on $U \times G$.
Thus $u_1,..., u_n$ give a trivialization of  $F_x^*T_Y$. \vskip
.2cm \noindent Let $Ram(f)'$ be the complex analytic subset of
$Ram(f)$ defined by $$Ram(f)' := \{ x \in Ram(f), g \cdot x \in
Ram(f) \mbox{ for each } g \in G \}.$$ We have established that
$F_x^*T_Y$ is trivial unless $x \in Ram(f)'$. If $Ram(f)'$
contains no component of codimension 1, then we are done by
setting $E$ to be the union of the singular locus of $X$ and
$Ram(f)'$. Assume that there exists a component $R$ of codimension
1 in $Ram(f)'$. $R$ is preserved by the $G$-action, namely,
$\Phi(R \times G) = R$. Let $B = f(R)= F(R \times G)$ be the
reduced irreducible divisor on $Y$. \vskip .2cm \noindent Let $x$
be a smooth point of $R$ with the following conditions.
\begin{enumerate}
\item$f$ has rank $n-1$ at  $g
\cdot x$ for some $g \in G$.
\item $R$ is the only irreducible component of $f^{-1}(B)$ which
contains the $G$-orbit $G \cdot x$.
\end{enumerate}
Note that the set of points of $R$ which does not satisfy (1) or
(2) must be of codimension $\geq 2$ and is preserved by $G$. Thus
if we show $F_x^*T_Y$ is trivial for $x$ as above, the proof of
Proposition 3.2 is complete. \\
Consider the exact sequence of sheaves of differentials
 $$ 0 \rightarrow \Omega^1_{B/Y}
\longrightarrow \Omega^1_Y|_B \longrightarrow \Omega^1_B
\rightarrow 0.$$ Since $B$ is a reduced divisor in the complex
manifold $Y$, $ \Omega^1_{B/Y}$ is isomorphic to the invertible
sheaf $ {\sO}_B(-B).$ The condition (2) guarantees that we can
write on $ U  \times G$,
$$F^*(-B) = -m(R \times G) - H$$ as  divisors
where $m$ is a positive integer and $H$ is an effective divisor
such that $\{x\} \times G$ is not contained in the support of $H$.
Since $ {\sO}(R \times G)$ is a trivial invertible sheaf along
$\{x\} \times G$, $$F_x^* {\sO}_B(-B) = F^*_x(-H).$$ This implies
$$\dim H^0(  G, F_x^*\Omega^1_{B/Y}) = \dim H^0(\{x\} \times G,{\sO}(F^*(-B))) \leq 1.$$
 Pulling back the sequence of differentials
 by $F_x$, we get the exact
sequence $$ 0 \rightarrow H^0(G,F_x^*\Omega^1_{B/Y} )
\longrightarrow H^0(G, F_x^* \Omega^1_Y) \longrightarrow H^0( G,
F_x^* \Omega^1_B).$$
 Note that when $z$ is a general point
 of $U$, $\dim H^0( G, F_z^* \Omega^1_Y) = n$ because we established
previously that  $F_z^*T_Y$
 is trivial. It follows from the upper-semi-continuity,  $$\dim H^0( G,
 F_x^* \Omega^1_B) \geq n-1.$$
Taking duals of the exact sequence of differentials,
 we get $$0 \rightarrow (\Omega^1_B)^*
\longrightarrow T_Y|_B \longrightarrow (\Omega^1_{B/Y})^*. $$ Here
$ (\Omega^1_B)^*$ is the subsheaf of the tangent sheaf of $Y$
consisting of vector fields tangent to $B$. \vskip .2cm \noindent
Now let us go back to our construction of sections $u_1, \ldots,
u_{n-1}$ of $F^*T_Y$ on $U \times G$. We can choose the
neighborhood $U$ and vector fields $v_1, \ldots, v_n$ on $U$ such
that all points of $U \cap R$ satisfy (i) and (ii), and  $v_1,
\ldots, v_{n-1}$ are tangent to $U \cap R$. Then the sections
$u_1, \ldots, u_{n-1}$ of $F^*T_Y$ on $U \times G$ will be
sections of the subsheaf $F^* (\Omega^1_B)^*$. By (i),  $B$ will
be smooth at $y:= f(g \cdot x)$ and  the differential $df$ at
$g\cdot x$ will send $T_{g \cdot x}(R)$ isomorphically into the
subspace $T_y(B)$ of  $T_y(Y)$. Thus $u_1, \ldots, u_{n-1}$ will
be pointwise independent at $(x, g)$ spanning the fiber $T_y(B)$.
Therefore the sections $u_1, \ldots, u_{n-1} \in H^0(G, F_x^*
(\Omega^1_B)^*)$ span the fiber at one point. \vskip .2cm
\noindent Together with $\dim H^0(G, F_x^* \Omega^1_B) \geq n-1$,
this implies that $u_1, \ldots, u_{n-1}$ as sections of $F_x^*
T_Y$
 will be pointwise independent at every point of $G$.
Thus they span a trivial subbundle $W_x$  of rank $n-1$ in $F_x^* T_Y$.
But we have seen that when $z$ is a general point of $U$,
$u_1, \ldots, u_{n-1}$ span a trivial subbundle $W_z$ of rank $n-1$
in $F^*_zT_Y$ such that the quotient is a trivial line bundle.
It follows that the quotient $F_x^*T_Y/W_x$ is also a trivial line
bundle. Thus we can write $F_z^*T_Y$ as an extension of the trivial
line bundle by a trivial vector bundle $W_z$ for each $z \in U$.
For general $z \in U$, we know that this extension is trivial.
Thus the extension is trivial for $x$, too. This completes the proof.
\end{proof}

Coming back to the set-up of the beginning of this section,
we find an unramified covering $h:Z \to Y$ and a ramified covering $g: X \to Z$
such that $f = h  \circ g$ such that
$$ df(v') \in g^*H^0(Z,T_Z) = g^*H^0(Z,h^*T_Y)$$
proving Theorem 1.3.

\section{Proof of the Main Result}
\setcounter{lemma}{0}

To prove Theorem 1.2, after Stein factorization, we may assume
that $f$ is finite. In particular, $X$ belongs to the Fujiki's
class ${\mathcal C}$ (e.g. [CP94, Proposition 3.17]). We fix in this
section a normal compact complex space $X$ of class ${\mathcal
C}$, a compact K\"ahler manifold $Y$ and  a finite surjective map
$$f: X \to Y.$$ This implies that each component of $\Hol(X,Y)$ is a Zariski
open subset in a compact complex space.  We always may assume that
$X$ and $Y$ are non-algebraic.

\begin{construction} {\rm  By [Ca81] there exists an almost holomorphic map
$$ q_X: X \dasharrow Q_X$$
to a compact K\"ahler manifold $Q_X$ with the following property: given two very general points $x,y \in X$,
we have $q_X(x) = q_X(y) $ if and only if $x$ and $y$ can be joined by a chain of irreducible curves, all of the
components of the chain belonging to families of curves which cover $X.$ \\
Of course, if $X$ is Moishezon, then $Q_X$ is a single point;
conversely by [Ca81], if $Q_X$ is a point, then $X$ is Moishezon.
If there is no covering family of curves in $X$; then $Q_X = X$
(up to birational transformation). Notice that $q_X$ is ``in
general'' different from the algebraic reduction of $X,$ the fiber
of which not necessarily
being algebraic. \\
In the same way we obtain a map $q_Y: Y \dasharrow Q_Y.$ \vskip
.2cm \noindent Clearly $f$ maps a general $q_X$-fiber to a general
$q_Y$-fiber, so that we obtain a commutative diagram
 \[\xymatrix{X \ar[d]^{q_X} \ar[r]^{f} \ar[d] & Y \ar[d]^{q_Y}\\
             Q_X \ar[r]^{\tilde f} & Q_Y}.\]

By
$$ \Hol_f^v(X,Y) \subset \Hol_f(X,Y)$$
we will denote the space of ``vertical'' deformations of $f,$
which means that we deform $f$ with $\tilde f$ fixed. The tangent
space to $ \Hol_f^v(X,Y)$ at $[f]$ is the subspace
$$ H^0(X,f^*T_{Y/Q_Y}) \subset H^0(X,f^*T_Y)$$ consisting of sections of $f^*T_Y$ which are tangent
 to the general fiber of $q_Y$ (and therefore tangent to all ``fibers'' of $q_Y.$).
}
\end{construction}

\begin{proposition} $\Hol_f^v(X,Y)$ is compact.
\end{proposition}

\begin{proof} We introduce the shorthand $H = \Hol_f^v(X,Y)$. Then $H$ can be considered as subset of the cycle space
$\sC(X \times Y).$ Since the irreducible components of $\sC(X \times Y)$ are compact (see [CP94] for a discussion and further references),
it suffices to prove that $H$ is closed in $\sC(X \times Y).$ \\
First we observe that we may assume $X$ smooth and K\"ahler. In
fact,  take a bimeromorphic holomorphic map $\pi: \hat X \to X$
such that $\hat X$ is smooth and K\"ahler. Let $\hat f = f \circ
\pi$ and $\hat H$ be the space of vertical deformations of $\hat
f.$
Then $H \subset \hat H$ in a natural way, and it is clear that $H$ is closed in $\sC(X \times Y)$ once we have proved that $\hat H$ is closed
in $\sC(\hat X \times Y)$. \\
To prove that $H$ is closed, we consider a family $(f_t)_{t  \in
\Delta^*}$ in $H$ over the punctured unit disc $\Delta$ such that
the
graphs $G_t \subset X \times Y$ converge to a cycle $G_0.$ We need to prove that $G_0$ is the graph of a holomorphic map $f_0: X \to Y.$  \\
Let us fix K\"ahler forms $\omega_X$ and $\omega_Y$ on $X$ and $Y.$ Let $p_X: X \times Y \to X$ and $p_Y: X \times Y \to Y$ be the projections. \\
We show first that $G_0$ is irreducible. In fact, $G_0$ has a unique component with multiplicity $1$, say $G^*,$ which maps onto $X$. Moreover
$G^* \to X$ has degree $1.$ Both statement follow by integrating $p_X^*(\omega_X)$ over $G_0$. \\
Now $G^*$ defines at least a meromorphic map $f_0: X \dasharrow Y.$ Since $\kappa (X_b) \geq 0$ for general $b  \in Q_X,$ the space
$$\Hol_{f \vert X_b}(X_b,Y_{\tilde f(b)})$$
is compact by [HKP05] and therefore $f_0$ is holomorphic near $X_b,$ i.e. $G^* \cap X_b$ is the graph of $f_0 \vert X_b.$ Since
$f_0 \vert X_b$ has degree $d = \deg f_t,$ we conclude that $G^* \to Y$ has degree $d.$
Introducing the K\"ahler form
$$ \omega = p_X^*(\omega_X) + p_Y^*(\omega_Y), $$
we notice that
$$ G^* \cdot \omega^n =  G_t \cdot \omega^n, $$
on the other hand $G_t \cdot \omega^n = G_0 \cdot \omega^n, $ so that
$$ G^* = G_0, $$
and $G_0$ is irreducible. \\
We consider the family $\sG = (G_t) \to \Delta$ with projection $p: \sG \to X \times \Delta. $
$ \sG$ is an irreducible reduced complex space and $\deg p = 1.$ Therefore $p$ has connected fibers, i.e.
$p$ is bimeromorphic. If $p$ were not biholomorphic, then the purity-of-branch theorem (recall that $X$ is smooth) exhibits a
proper subspace $B \subset X \times \Delta,$ necessarily contained in $X \times \{0\},$ such that $D = p^{-1}(B)$ has codimension 1
in $\sG,$ i.e., $\dim D = n.$ Hence $D$ is an irreducible component of $p_{\Delta}^{-1}(0).$ This contradicts the irreducibility of $G_0.$ \\
Hence $p$ is biholomorphic, and so does $G_0 \to X.$ Therefore $f_0$ is holomorphic and of course $f_0 \in \Hol_f^v(X,Y),$ which completes
the proof.
\end{proof}

\begin{proposition} Let $f: X \to Y$ be a surjective holomorphic map to a compact K\"ahler manifold $Y$ of non-negative Kodaira dimension.
Suppose there exists a positive dimensional compact subvariety in $\Hol_f(X,Y).$ Then $f$ has an unramified factorization $\beta: Z \to Y$
with $\dim {\rm Aut}(Z) > 0.$
\end{proposition}

\begin{proof}  By Proposition 2.5, we know that
there exists a (not necessarily unramified) factorization $h: Y'
\rightarrow Y$ with $\dim {\rm Aut}(Y') > 0$. Since $\kappa (Y) >\geq 0,$ the identity component ${\rm Aut}^0(Y')$ is a torus $G.$
Applying Theorem 1.3 to $h$, we get an unramified factorization
$\beta:Z \rightarrow Y$ with an effective $G$-action on $Z$. This
gives an unramified factorization of $f$ with $\dim {\rm Aut}(Z) >0$.
\end{proof}

\begin{Proof of Theorem 1.2} {\rm To prove Theorem 1.2, we
we may assume by Proposition 2.3, by Proposition 2.4 and by
induction on $\dim Y$, that $f$ has no unramified factorization
$\beta: Z \rightarrow Y$ with $\dim {\rm Aut}(Z) >0$. Then by
Proposition 4.3, it follows that $\Hol_f(X, Y)$ has no positive
dimensional compact subvariety. Consequently by Proposition 4.2,
$\dim Hol_f^v (X, Y) = 0$. \vskip .2cm \noindent Now consider
closed irreducible analytic set (not necessarily compact)
$$ Z \subset \Hol_f(X,Y) $$
(e.g., $Z = \Hol_f(X,Y)$) and argue as in [No92, Section 3].  For
$g \in Z$ consider the ramification locus $B_g \subset Y.$ If $B_g
\subset Y$ (resp. $g^{-1}(B_g) \subset X$) moves, consider a
general point $y \in Y$ (resp. $x \in X$) and form the
1-codimensional subvariety $Z(x) \subset Z$ consisting of those
$g$ with $y \in B_g$ (resp. $x  \in g^{-1}(B_g)).$ Repeating this
process we obtain the following alternative.
\begin{itemize}
\item $\dim Z = 1;$ \item there exists a positive-dimensional $Z$
such that $B_g$ and $g^{-1}(B_g)$ are independent of $g$ for all
$g \in Z.$
\end{itemize}

In the first case we consider for fixed general $x$ the closed
curve $C = \{g(x) \vert g \in Z \}.$ Since $\Hol_f(X,Y)$ is
Zariski-open in the cycle space (see e.g. [Fu78], [Li78], at least
in the case $f$ is an automorphism), we can take closure and
obtain a compact curve $\overline C$ parametrizing cycle such that
the general cycle is the graph of a deformation of $f.$ Now take
$x \in X$ general and consider the closure of the curve
$$C_x = \{g(x) \vert  g \in C\}. $$
Thus we obtain a covering family of compact curves in $X$ and
since the deformations of $f$ are not vertical, these curves are
not contained in fibers of $q_X.$ This contradicts the
construction of $q_X.$ \vskip .2cm \noindent In the second case we
may apply Proposition 2.6 and $Z$ is (isomorphic to) a closed
subvariety of the automorphism group of $X$, which is a torus, a
contradiction. }
\end{Proof of Theorem 1.2}

\bigskip
{\bf Acknowledgment.} We would like to thank Junjiro Noguchi for
the references [Ho87] and [No92].

\end{document}